\documentclass[12pt,leqno,oneside]{amsart}
\usepackage{mathrsfs,dsfont}
\usepackage{amsmath,amstext,amsthm,amssymb,bbm,faktor}
\usepackage{charter}
\usepackage{typearea}
\usepackage{hyperref}
\usepackage{color}
\usepackage{nccmath}
\usepackage{float}
\usepackage{tikz}
\usetikzlibrary{patterns}
\usepackage{graphicx}
\usepackage{hyperref}
\usepackage{cleveref}

\usepackage[width=6.4in,height=8.5in]
{geometry}

\def\bt{\begin{thm}}
	\def\et{\end{thm}}
\def\bl{\begin{lem}}
	\def\el{\end{lem}}
\def\bd{\begin{defn}}
	\def\ed{\end{defn}}
\def\bc{\begin{cor}}
	\def\ec{\end{cor}}
\def\bp{\begin{proof}}
	\def\ep{\end{proof}}
\def\br{\begin{rem}}
	\def\er{\end{rem}}
\def\brs{\begin{rems}}
	\def\ers{\end{rems}}

\newtheorem{thm}{Theorem}[section]
\newtheorem{prop}[thm]{Proposition}
\newtheorem{lem}[thm]{Lemma}
\newtheorem{defn}[thm]{Definition}

\newtheorem{rem}[thm]{Remark}
\newtheorem{rems}[thm]{Remarks}
\newtheorem{cor}[thm]{Corollary}

\numberwithin{equation}{section}

\title{The Smirnov Class for Sub-Bergman Hilbert Spaces}

\author{Sibel \c{S}ahin} 

\address{Mathematics Department, Mimar Sinan Fine Arts University, \.{I}stanbul, Turkey}
\email{sibel.sahin@msgsu.edu.tr}

\begin{document}
	
	\begin{abstract}
		In this work we consider the Smirnov classes for sub-Bergman spaces. First we point out some observations about the Smirnov property of sub-Bergman space $H^\alpha(\varphi)$ and its relation to the defining $H^{\infty}_{1}$ function $\varphi$. The first main result of the work deals with the range of the defect operator over Smirnov-sub-Bergman class whereas in the last part we show that contrary to classical Bergman spaces, Smirnov-sub-Bergman classes have non-tangential boundary values almost everywhere on the unit circle.
		
	\end{abstract}
	
	\maketitle
	%\tableofcontents

	\section{Introduction}
	
	  In the classical theory of harmonic or analytic Hardy spaces it is known that any harmonic or analytic Hardy function ($h^p$ or $H^p$ respectively) has radial limits almost everywhere in the boundary for $p\geq 1$. This fact happens to be true also for analytic functions of the unit disc belonging to $H^p$ classes for $p<1$ and this type of generalization in fact led to the investigation of a wider class of analytic functions namely the Nevanlinna class $N$. A famous theorem of F. and R. Nevanlinna states that an analytic function of the unit disc belongs to $N$ if and only if it is the quotient of two bounded analytic functions (\cite{Dur}, Theorem 2.1, pg:16). Furthermore, due to a well-known theorem by Smirnov any analytic function $f\in H^2(\mathbb{D})$ can be written in the form $f=\varphi/\psi$ where $\varphi$ and $\psi$ are both multipliers (i.e they are both from $H^\infty(\mathbb{D})$) and $\psi$ is an outer function hence $\psi H^2$ is dense in $H^2(\mathbb{D})$ (i.e $\psi$ is a cyclic multiplier). Inspired by this, in \cite{AHMR}, Aleman et al. defined the Smirnov class $\mathcal{N}^+$ of a general reproducing kernel Hilbert space $\mathcal{H}$ as the set of quotients $\varphi/\psi$ where $\varphi$ and $\psi$ are multipliers of $\mathcal{H}$ and $\psi$ is a cyclic multiplier of $\mathcal{H}$. It is clear that the underlying reproducing kernel Hilbert space changes the behavior of the multipliers and as a consequence the structure of Smirnov class changes drastically so let us mention some of these different reproducing kernel Hilbert spaces of analytic functions in the unit disc. 
	
	  Let $T$ be a bounded linear operator on a Hilbert space $H$ and suppose that $T$ is a contraction then we have an inner product space as the range of so called defect operator $(I-TT^*)^{1/2}$ and it is denoted by $H(T)$. The inner product in consideration is defined as follows:
	
	$$
	\left\langle (I-TT^*)^{1/2}x,(I-TT^*)^{1/2}y \right\rangle_{H(T)}=\langle x,y \rangle_{H}
	$$      	
	for $x,y\in H\ominus Ker(I-TT^*)^{1/2}$.

	When the underlying Hilbert space is the Hardy space $H^2(\mathbb{D})$ and the aforementioned contraction is the Toeplitz operator $T_b$ for some $b\in H_{1}^{\infty}(\mathbb{D})$ then $H(T)=H(b)$ is called de Branges- Rovnyak space. For an extensive study of these spaces one may refer to (\cite{Fri1}, \cite{Fri2}). When the underlying Hilbert space is a generalized Bergman space $A^2_\alpha$ and the contraction is the analytic Toeplitz operator $T^\alpha_\varphi:A^2_\alpha\longrightarrow A^2_\alpha$ with symbol $\varphi$ belonging to the unit ball of multiplier algebra $\mathcal{M}(A^2_alpha)$ of $A^2_\alpha$ then $H(T)=H^\alpha(\varphi)$ is called a sub-Bergman space. For a comprehensive study of sub-Bergman spaces one may see (\cite{Zhu1}, \cite{Zhu2}, \cite{Zhu3}). Since Aleman's generalization of Smirnov classes can be applied to any reproducing kernel Hilbert space where one knows the behavior of multiplier algebra, one can also consider the Smirnov classes of de-Branges-Rovnyak spaces (\cite{FHRT}) or as we did in this work the Smirnov classes of sub-Bergman spaces.
	
	A space $H$ is said to be satisfying the Smirnov property if $H\subset \mathcal{N}^{+}(H)$ which is in fact the case for the classical Hardy space $H^2(\mathbb{D})$ and this result can be shown using the existence of inner-outer factorization and radial limits. For the de Branges-Rovnyak case $H=H(b)$, it is shown in \cite{FHRT} that
	
	\bt (\cite{FHRT}, Theorem 1.1): If $b\in H^\infty_1$ is rational but not inner then $H(b)\subset \mathcal{N}^{+}(H(b))$.
	\et
	
	In \cite{AHMR}, Aleman et al. showed that actually the Smirnov property of the reproducing kernel space $H$ is highly connected to the kernel itself but before mentioning this remarkable result let us mention what a complete Nevenlinna-Pick kernel is:
	
	\bd
	Let $H$ be a reproducing kernel Hilbert space on a set $X$. Then $H$ is called to be a complete Pick space if for any $r\in \mathbb{N}$; $z_1,\dots,z_n\in X$ and the matrices $W_1,\dots, W_n\in M_r(\mathbb{C})$, positivity of the $nr\times nr$ matrix 
	$$
	[K(z_i,z_j)(I_{\mathbb{C}^r}-W_iW^*_j)]_{i,j=1}^{n}
	$$ 
	gives the existence of $\Phi\in M_r(Mult(H))$ of norm at most $1$ such that
	$$
	\Phi(z_i)=W_i,~~i=1,\dots, n.
	$$
	In this case the kernel $K(z,w)$ is called a complete Nevanlinna-Pick kernel.
	\ed      

	In \cite{AHMR}, Aleman et al. showed that if the reproducing kernel of the Hilbert space $H$ satisfies the complete Nevanlinna-Pick property along with normalization i.e $K(z,z_0)=1$ for some $z_0\in \mathbb{D}$ and for all $z\in\mathbb{D}$ then $H\subset \mathcal{N}^{+}(H)$.
	
	In this paper, we will consider the Smirnov class for the sub-Bergman Hilbert spaces. Details of sub-Bergman Hilbert spaces, for both classical and weighted Bergman cases, together with a detailed study of their Nevanlinna-Pick property can be found in (\cite{Zhu1},\cite{Zhu2}, \cite{Zhu3}).  

    \section{Sub-Bergman Spaces and Smirnov Class}

	Before studying the general Smirnov classes let us first mention the classical Nevanlinna spaces (also known as Smirnov class or functions of bounded characteristic):
	
	\bd\label{originaldef}
	A holomorphic function of the unit disc $\mathbb{D}$ is said to be of class $\mathcal{N}$ if the integrals 
	$$
	\int_{0}^{2\pi} \log^{+}|f(re^{i\theta})|d\theta
	$$
	are bounded for $r<1$. 
	\ed
	
	By (\cite{Dur}, Theorem 2.9, pg:25) we know that any $f\not\equiv 0$ of class $\mathcal{N}$ can be factorized like $f=B(z)[S_1(z)/S_2(z)]F(z)$ where $B(z)$ is a Blaschke product, $S_1$ and $S_2$ are singular inner functions and $F$ is an outer function. Moreover, the class $\mathcal{N}^{+}$ is the set of all functions $f\in \mathcal{N}$ for which $S_2\equiv 1$. The class $\mathcal{N}^{+}$ is the natural limit of Hardy classes $H^p$ as $p\rightarrow 0$ with the proper inclusions $H^p\subset \mathcal{N}^{+}\subset \mathcal{N}$.
	
	For a general Hilbert function space $H$ on a set $X$, the general Smirnov class is defined as follows:
	
	\bd
	Let $H$ be a reproducing kernel Hilbert space of analytic functions on $\mathbb{D}$, then the multiplier algebra is defined as
	$$
	Mult(H)=\{\varphi\in\mathcal{O}(\mathbb{D}):~~\varphi H\subseteq H\}
	$$
	and the Smirnov class of $H$ is defined as
	$$
	\mathcal{N}^{+}(H)=\left\{\frac{\varphi}{\psi}:~~\varphi,\psi\in Mult(H)~\text{and}~\overline{\psi H}=H\right\}
	$$
	where the closure is with respect to the norm of $H$.
	\ed

	Since we are interested in the Smirnov class of sub-Bergman spaces in this work, let us first give a brief introduction to Bergman and sub-Bergman spaces:
	
	Let $A^{2}_{\alpha}$ denote the space of analytic functions $f$ on $\mathbb{D}$ such that 
	$$
	\int_{\mathbb{D}}|f(z)|^2(1-|z|^2)^\alpha dA(z)<\infty,~~\alpha> -1.
	$$

	Then $A^{2}_{\alpha}$ (Generalized Bergman Space) is a reproducing kernel Hilbert space with the following inner product
	$$
	\langle f,g\rangle=\int_{\mathbb{D}}f(z)\overline{g(z)}(1-|z|^2)^\alpha dA(z)
	$$
	and the reproducing kernel 
	$$
	K(z,w)=K^{\alpha}_{w}(z)=\dfrac{1}{(1-z\overline{w})^{2+\alpha}}.
	$$
	
	It is known that for $\alpha \geq 1$, $Mult(A^{2}_{\alpha})=H^\infty(\mathbb{D})$. For $\varphi\in Mult(A^{2}_{\alpha})$, Toeplitz operator $T_\varphi:A^{2}_{\alpha}\longrightarrow A^{2}_{\alpha}$ is a contraction and the range of the defect operator $D_{\varphi}^{\alpha}=(I-T_\varphi T_{\varphi}^{*})^{1/2}$ is called \textit{the sub-Bergman space} and denoted by $H^{\alpha}(\varphi)$. In \cite{Zhu3}, it has been shown that $H^{\alpha}(\varphi)$ is itself a reproducing kernel Hilbert space with the kernel
	$$
	K^{\alpha,\varphi}(z,w)=K^{\alpha,\varphi}_{w}(z)=\dfrac{1-\varphi(z)\overline{\varphi(w)}}{(1-z\overline{w})^{2+\alpha}}.
	$$
	
	Throughout this work we will be interested in the class $\mathcal{N}^{+}(H^{\alpha}(\varphi))$. We have already mentioned the result of Aleman et al. \cite{AHMR} relating the complete Nevanlinna-Pick kernel to the Smirnov property of the given space. A recent study of Luo and Zhu (\cite{Zhu3}) gives the complete characterization for the complete Nevanlinna-Pick kernels for sub-Bergman spaces $H^{\alpha}(\varphi)$ as follows:
	
	\bt (\cite{Zhu3}, Theorem 5)
	Suppose $\varphi\in H_{1}^{\infty}$ and $-1<\alpha\leq 0$. Then the reproducing kernel of $H^{\alpha}(\varphi)$ is a complete Nevanlinna-Pick kernel if and only if $\varphi$ is a M\"{o}bius map.
	\et
	
	Hence combining this with the result of Aleman et al. we obtain the following:
	
	\begin{prop}
	Let $-1<\alpha\leq 0$ and $\varphi$ be a M\"{o}bius map then $H^{\alpha}(\varphi)\subset \mathcal{N}^{+}(H^{\alpha}(\varphi))$.
	\end{prop}
	
	In \cite{Zhu2} and \cite{Zhu3}, apart from the results related to Nevanlinna-Pick property the range of the defect operators were studied for both classical and generalized Bergman spaces. In the classical case, defect operator $(I-T_BT_{B^*})^{1/2}$ maps the Bergman space $A^2(\mathbb{D})$ to the Hardy space $H^2(\mathbb{D})$ (\cite{Zhu2}, Theorem A, pg:328) and in the generalized case the space $A^{2}_{\alpha}$, $\alpha>-1$, is mapped to $A^{2}_{\alpha-1}$ by the defect operator $D^{\alpha}_{\varphi}$ where $\varphi$ is a finite Blaschke product (\cite{Zhu3}, Theorem 11, pg:13).
	
	Now in the following first main result of this study, we will consider the range of the defect operator on the Smirnov-sub-Bergman classes $\mathcal{N}^{+}(H^{\alpha}(\varphi))$:
	
	\bt (Main Theorem 1)
	Let $\alpha>-1$ and $\varphi$ be a M\"{o}bius map. Then the defect operator $D^{\alpha}_{\varphi}=(I-T_\varphi T_{\varphi^*})^{1/2}$ maps the Smirnov-sub-Bergman class $\mathcal{N}^{+}(H^{\alpha}(\varphi))$ to $\mathcal{N}^{+}(H^{\alpha-1}(\varphi))$.
	
	In particular, the Smirnov-Bergman class $\mathcal{N}^{+}(A^2(\mathbb{D}))$ is mapped to classical Nevanlinna (Smirnov) class $\mathcal{N}^+$ ($\ref{originaldef}$).
	\et
	
	In order to prove the main theorem 1, we need the following lemma:
	
	\bl (\cite{Zhu2}, Proposition 3)
	Let $B$ be a finite Blaschke product with distinct zeros. If $f\in A^2(\mathbb{D})$ and $F$ is any anti-derivative of $f$ in $\mathbb{D}$, then
	$$
	T_{\overline{B}}f(z)=\dfrac{f(z)}{B(z)}-\dfrac{F(z)B'(z)}{(B(z))^2}+\sum_{k=1}^{N}\dfrac{F(a_k)}{B'(a_k)(z-a_k)^2}
	$$
	where $a_1,\dots,a_N$ are the zeros of $B$ in $\mathbb{D}$. 
	
	\el
	
	\bp (of Main Theorem 1)
	Let $f\in \mathcal{N}^{+}(H^{\alpha}(\varphi))$ and $\varphi$ be a M\"{o}bius map. Then $f=\dfrac{\gamma}{\psi}$ for some $\gamma$ and $\psi$ such that $\gamma,\psi\in Mult(H^{\alpha}(\varphi)$ and $\overline{\psi H^{\alpha}(\varphi)}=H^{\alpha}(\varphi)$ and without loss of generality $\varphi(z)=\xi\dfrac{a-z}{1-\overline{a}z}$ with $\xi\in\partial\mathbb{D},~a\in\mathbb{D}$.
	
	Now by previous lemma,
	$$
	T_{\overline{\varphi}}\left(\dfrac{\gamma}{\psi}\right)=	T_{\varphi^*}\left(\dfrac{\gamma}{\psi}\right)=\dfrac{\gamma}{\psi\varphi}-\dfrac{F\varphi '}{\varphi^2}+\dfrac{F(a)}{\varphi'(a)(z-a)^2}
	$$
	where $F$ is any anti-derivative of $f$ in $\mathbb{D}$ and 
	$$
	T_{\varphi}(T_{\overline{\varphi}})=T_{\varphi}T^{*}_{\varphi}\left(\dfrac{\gamma}{\psi}\right)=\dfrac{\gamma}{\psi}-\dfrac{F\varphi'}{\varphi}+\dfrac{F(a)\varphi}{\varphi'(a)(z-a)^2}
	$$ 
	since $T_\varphi(g)=P(\varphi g)$ but we have analytic $\varphi$ so it is just multiplication operator. 
	
	Hence 
	\begin{equation*}
		(I-T_{\varphi}T^{*}_{\varphi})\left(\dfrac{\gamma}{\psi}\right)=\left(\dfrac{F\varphi '}{\varphi}-\dfrac{F(a)\varphi}{\varphi'(a)(z-a)^2}\right)
	\end{equation*}
	\begin{equation}\label{eq:star}
		=\dfrac{F\varphi'\varphi'(a)(z-a)^2-F(a)\varphi^2}{\varphi\varphi'(a)(z-a)^2}
	\end{equation}
	Now from (\cite{Zhu3}, Theorem 11), we know that in this setting $H^{\alpha-1}(\varphi)=A^{2}_{\alpha-2}$ so in order to show that $(\ref{eq:star})$ belongs to $\mathcal{N}^+(H^{\alpha-1}(\varphi))$ we need to see that
	\begin{equation}\label{eq:doublestar}
	\varphi'(a)(z-a)^2\varphi A^{2}_{\alpha-2}(\mathbb{D})= A^{2}_{\alpha-2}(\mathbb{D})
	\end{equation} 
	is satisfied. For this observe that for the reproducing kernel $K_{w}^{(\alpha-1),\varphi}$ of $H^{\alpha-1}(\varphi)$ we have 
	$$
	\varphi'(a)(z-a)^2\varphi K_{w}^{(\alpha-1),\varphi}(z)
	$$
	$$
	=\dfrac{1-|a|^2}{(1-\overline{a}z)(1-a\overline{w})}\dfrac{1}{(1-z\overline{w})^\alpha}\dfrac{a-z}{1-\overline{a}z}(z-a)^2
	$$
	which is just a re-scaled version of $K^{\alpha-2}_{w}(z)$, hence we have $(\ref{eq:doublestar})$. Therefore the defect operator maps $\mathcal{N}^{+}(H^{\alpha}(\varphi))$ to the space $\mathcal{N}^{+}(H^{\alpha-1}(\varphi))$. 
	\ep
	
	As a direct consequence of this theorem we have the following slight generalization:
	
	\bc
	Let $\alpha>-1$ and $\varphi$ be a finite Blaschke product with distinct zeros then the defect operator $D^{\alpha}_{\varphi}=(I-T_\varphi T_{\varphi^*})^{1/2}$ maps the Smirnov sub-Bergman class $\mathcal{N}^{+}(H^{\alpha}(\varphi))$ to $\mathcal{N}^{+}(H^{\alpha-1}(\varphi))$. 
	\ec
	
	As it was mentioned at the beginning the main motivation behind defining the Smirnov (or Nevanlinna) classes was to understand the boundary behavior of the holomorphic functions of the unit disc belonging to less regular Hardy spaces and when it comes to the Bergman spaces as it was pointed out by Aleman et al. (\cite{AHMR}, pg:229), the classical Bergman space $A^{2}_{\alpha}(\mathbb{D})$ contains some functions that do not have radial boundary limits on $\mathbb{T}=\partial\mathbb{D}$. However, contrary to classical Bergman spaces, some Smirnov sub-Bergman classes have non-tangential boundary values almost everywhere on $\mathbb{T}$:
	
	\bt (Main Theorem 2)
	Let $f\in \mathcal{N}^{+}(H^{\alpha}(\varphi)) $ for some $\varphi\in Mult (A^{2}_{\alpha})(\mathbb{D})$. Then $f$ has non-tangential boundary values $f(\xi)$ for almost all $\xi\in\mathbb{T}$.
	\et  
	\bp
	Suppose that $f\not\equiv 0$ is of the form $f=\dfrac{\gamma}{\psi}$ where $\gamma$ and $\psi$ are bounded holomorphic functions of $\mathbb{D}$. Then by (\cite{Dur}, Theorem 1.3, pg:6) both $\gamma$ and $\psi$ have non-tangential limits $\gamma(\xi)$ and $\psi(\xi)$ for almost all $\xi\in\mathbb{T}$. In particular, $\psi(e^{i\theta})$ cannot be zero on a set of positive measure since it is a cyclic multiplier therefore the non-tangential boundary value $f(\xi)$ exists almost everywhere on $\mathbb{T}$.
	\ep
	
\end{document}